\documentclass[12pt]{article}
\usepackage{amsmath}
\usepackage{amsfonts}
\usepackage{amssymb}

\newtheorem{ex}{Example}{\bf}{\rm}{\rm}

\def\Real{\mathbb{R}}

\def\so{\mathfrak{so}}
\def\simil{\mathfrak{sim}}
\def\Simil{{\rm Sim }}

\newcommand{\tx}{\tilde{x}}
\renewcommand{\d}{{\rm d}}

\def\Ric{\mathop{{\rm Ric}}\nolimits}
\def\p{\partial}
\def\t{\tilde}

\def\CQG{Class. Quantum Grav.}

\def\be{\begin{equation}}
\def\ee{\end{equation}}

\begin{document}

\title{Examples of Einstein spacetimes with recurrent null
vector fields}
\author{Anton S. Galaev}

%\address{
%Department of Mathematics and Statistics, Faculty of Science,
%Masaryk University in Brno, Kotl\'a\v rsk\'a~2, 611~37 Brno, Czech
%Republic } \ead{galaev@math.muni.cz}

\maketitle

\begin{abstract} The Einstein Equation on 4-dimensional  Lorentzian manifolds admitting
recurrent null vector fields is discussed. Several examples of a
special form are constructed. The holonomy algebras, Petrov types
and the Lie algebras of Killing vector fields of the obtained
metrics are found.
\end{abstract}

%\pacs{04.20.Jb, 04.40.Nr}
%\maketitle

\section{Introduction}

We study the Einstein Equation on spacetimes  admitting parallel
distributions of null lines and construct several special
examples.

Let $(M,g)$ be a 4-dimensional Lorentzian manifold admitting a
parallel distribution of null lines (i.e. $(M,g)$ is {\it a Walker
manifold}). This condition holds if and only if the holonomy group
of $(M,g)$ is contained in the subgroup $\Simil(2)\subset{\rm
O}(1,3)$ preserving a null line in the Minkowski space
$\Real^{1,3}$. Equivalently, $(M,g)$ admits a recurrent null
vector field in a neighborhood of each point, and locally there
exist coordinates $v,x^1,x^2,u$ (so called {\em Walker
coordinates})  such that the metric $g$ has the form
\begin{equation}\label{Walker} g=2\d v\d u+h+2A\d u+H\cdot (\d u)^2,\end{equation}
where $h=h_{ij}(x^1,x^2,u)\d x^i\d x^j$ is an $u$-dependent family
of Riemannian metrics,\\  $A=A_i(x^1,x^2, u)\ \d x^i$ is an
$u$-dependent family of one-forms, and $H$ is a local function on
$M$, \cite{Walker}. The vector field $\p_{v}=\frac{\p}{\p v}$ is
null and recurrent, and it defines the parallel distribution of
null lines. We will also write  $x^1=x$, $x^2=y$.

A Lorentzian manifold $(M,g)$ is called {\it an Einstein manifold}
if $g$ satisfies the equation
\begin{equation}\label{Einstein} \Ric=\Lambda g,\qquad \Lambda \in \Real,\end{equation}
where $\Ric$ is the Ricci tensor of the metric $g$, i.e. $\Ric_{a
b}=R^c_{\ a  c b}$,
 where $R$ is the curvature tensor of the metric $g$.
 The number $\Lambda \in \Real$ is called {\it the  cosmological
 constant}. If $\Lambda=0$, i.e. $\Ric=0$, than the manifold is called {\it Ricci-flat
 or vacuum Einstein.} We consider the case $\Lambda\neq 0$.

The Walker manifolds are of particular type of Kundt spaces
\cite{BCH1,CHPP}. Recently  G.W.~Gibbons and C.N.~Pope \cite{G-P}
considered the Einstein equation on Walker manifolds of arbitrary
dimension.

A special example of the metric \eqref{Walker} is the metric of a
pp-wave. It is given by $h=(dx)^2+(dy)^2$, $A=0$ and $H$
independent of $v$. If such metric is Einstein, then it is vacuum
Einstein, and this happens whenever $(\p^2_x+\p^2_y)H=0$. For the
plane waves ($H=A_{ij}(u)x^ix^j$), sometime  it is useful to
rewrite the metric in the Rosen coordinates, $g=2dvdu+h$, where
$h=C_{ij}(u)x^ix^j$ is a family of flat metrics, see e.g.
\cite{BFP02,Ivashchuk,Schimming}. The examples of the Einstein
metrics that we construct have a similar structure.

Vacuum Einstein Walker metrics in dimension 4 are found in
\cite{KG2}. After a proper change of coordinates they are given by
$h=(dx)^2+(dy)^2$, $A_2=0$, $H=-(\p_xA_1)v+H_0$, where $A_1$ is a
harmonic function and $H_0$ can be found from a Poisson equation.
In \cite{GLE} several examples of such metrics are rewritten in
new coordinates such that $A=0$ and $h$ is an $u$-family of flat
metrics on $\Real^2$.

 In \cite{Lewandowski}, all
4-dimensional Einstein Walker metrics with $\Lambda\neq 0$ are
described. After a proper choice of the coordinates, $h$ becomes
an independent of $u$ metric  of constant curvature. Next,
$A=Wdz+\bar W\bar z$, $W=i\p_z L,$ where $z=x+iy$,  $L$ is
$\Real$-valued function  given by
\begin{equation}\label{Lphi}L=2\text{\rm Re}\left(\phi\p_z(\ln
P_0)-\frac{1}{2}\p_z\phi\right), \quad
2P_0^2=\left(1+\frac{\Lambda}{|\Lambda|}
z\bar{z}\right)^2,\end{equation} where $\phi=\phi(z,u)$ is an
arbitrary function holomorphic in $z$ and smooth in $u$. Finally,
$H=\Lambda^2v+H_0$, where the function $H_0=H_0(z,\bar{z},u)$ can
be found in a similar way.

An example of the Einstein Walker metric with $\Lambda\neq 0$ is
given in \cite{GT}. After a change of the coordinates in
\cite{G-P}, this metric is given by $A=0$ and $h$ independent of
$u$. In \cite{CGHP}, the universality of the metric from \cite{GT}
and more generally of Einstein metrics with ${\rm
Sim}(2)$-holonomy is proved.

{\it The aim} of this paper is to construct examples of Einstein
Walker metrics with $\Lambda\neq 0$ such that $A=0$ and $h$
depends on $u$. The solutions from \cite{Lewandowski} are not
useful for constructing examples of such form, since "simple"
functions $\phi(z,u)$  define complicated forms $A$. First we find
Walker metrics with $h$ independent of $u$ and $A\neq 0$, then we
change the coordinates in such a way that $A=0$. The constructed
examples can be useful according to \cite{BCH1,CGHP,CHPP,G-P}.
Similar examples can be  constructed  in dimension 5, this
dimension is discussed e.g. in \cite{G-P,Grover}.

\section{Coordinate transformations and reduction of the Einstein
equation}\label{Seccoord}

The Einstein Equation \eqref{Einstein} for the Walker metric
\eqref{Walker} in arbitrary dimension is written down in
\cite{G-P}. It is shown that $h$ is a family of Einstein
Riemannian metrics with the cosmological constant $\Lambda$,
$H=\Lambda v^2+H_1v+H_0$, $\p_vH_1=\p_vH_0=0$, and there are three
additional equations.

The Walker coordinates are not defined canonically and any other
Walker coordinates $\tilde v,\tx^1, \tx^2,\tilde u$ such that
${\p}_{\t v}=\p_{v}$ are given by the  transformation
\be\label{trans} \tilde v=v+f(x^1,x^2,u),\quad \tilde x^{i}=\t
x^i(x^1,x^2,u),\quad \tilde u=u+c,\ee see \cite{Schimming,GLE}. If
the metric \eqref{Walker} is Einstein with $\Lambda\neq 0$, then
there exist Walker coordinates such that $A=0$ and $H_1=0$
\cite{GLE}. Then the Einstein Equation takes the form
\be\label{EEred} \Delta H_0+\frac{1}{2}h^{ij}\ddot{h}_{ij}=0,\quad
\nabla^j\dot h_{ij}=0,\quad h^{ij}\dot{h}_{ij}=0, \quad
\Ric_{ij}=\Lambda h_{ij},\ee where
$\dot{h}_{ij}=\partial_u{h_{ij}}$. Thus we get two equations on
the family of Riemannian Einstein metrics and the Poisson equation
for the function $H_0$. Of course, a solution gives the constant
family, $\dot{h}_{ij}=0$ and a harmonic $H_0$, but it is more
interesting to find solutions with $\dot{h}_{ij}\neq 0$. This will
be done below.

In order to find the coordinates as in \cite{GLE}, we start with
any Walker coordinates and consider two transformations. First
consider the transformation $v\mapsto v+\frac{1}{2\Lambda}H_1$,
after that $H_1=0$. Then consider the transformation
$x^i\mapsto\tilde x^i(x^1,x^2,u)$ such that the inverse
transformation satisfies the system of ordinary differential
equations $$\frac{\d x^i(u)}{d u}=W^i(x^1(u),x^2(u),u),$$ where
$W^i=- A_jh^{ij}$, with the initial conditions $x^i(u_0)=\tilde
x^i$. The solution can be then written as $x^i=x^i(\t x^1,\t x^2,
u)$.

  Note that in dimension 2 (and 3) any Einstein Riemannian metric
  has constant sectional  curvature, hence any such metrics with the
  same $\Lambda$ are locally isometric, and the coordinates can be
  chosen in such a way that $\p_u h=0$. As in \cite{Lewandowski}, it is not hard to show
  that if $\Lambda>0$, then we may choose $h=(\d x)^2+\sin^2 x\,(\d
  y)^2$, $H=\Lambda v^2+H_0$, and the Einstein Equation is
  reduced to \be\label{fL>0} \Delta_{S^2}f=-2f,\,
 \Delta_{S^2}H_0=2\Lambda\left(2 f^2-(\p_xf)^2+\frac{(\p_yf)^2}{\sin^2
 x}\right),\,
\Delta_{S^2}=\p_x^2+\frac{\p_y^2}{\sin^2 x}+\cot x\,\p_x. \ee
 The function $f$
defines the 1-form $A$: $A=-\frac{\p_yf}{\sin x\, }\d x+ \sin x\,
\p_xf\d y$. Similarly, if $\Lambda<0$, then we chose
$h=\frac{1}{-\Lambda\cdot x^2}\big((\d x)^2+(\d y)^2\big)$, and
get \be\label{f,L<0} \Delta_{L^2}f=2f,\, \Delta_{L^2}H_0=-4\Lambda
f^2-2\Lambda x^2((\p_xf)^2+(\p_yf)^2),\,
\Delta_{L^2}=x^2(\p_x^2+\p_y^2),\ee and $A=-\p_yf\d x+\p_xf\d y$.
 Thus in order to find partial solutions of \eqref{EEred} it is
 convenient to solve one of the above equations for $f$ and then, changing
 the coordinates, to get rid of $A$. Note that  {\it after such change
 $h$ does not depend on $u$ if and only if $A$ is a Killing form
 for $g$} \cite{G-P}. For $\Lambda>0$ this happens if and only if
 $f=c_1(u)\sin x\sin y+c_2(u)\sin x\cos y+c_3(u)\cos x$;
for $\Lambda<0$ this happens if and only if
$f=c_1(u)\frac{1}{x}+c_2(u)\frac{y}{x}+c_3(u)\frac{x^2+y^2}{x}$.
The functions $\phi(z,u)=c(u),$ $c(u)z$ and $c(u)z^2$ from
\eqref{Lphi} define a Killing form $A$, see \cite{GLE}, and  $A$
becomes complicated for other $\phi$.

Let now $g$ be a Walker Einstein metric \eqref{Walker} with $A=0$,
and let  $\tilde h(u)$ be a family of Riemannian metrics such that
$h(u_0)=\tilde h(u_0)$ for some $u_0\in\Real$. Then {\it  the
metric $g$ is isometric to the metric $\tilde g=2\d v\d u+\tilde
h+\tilde H(du)^2$ with some $\tilde H$ if and only if $h=\tilde h$
for all $u$.} Indeed, if the metrics are isometric, then one can
be taken to another using transformation \eqref{trans}. Since
$h(u_0)=\tilde h(u_0)$, we may assume that $\tilde
x^k(x^1,x^2,u_0)=x^k$.
 Note that $A_i=\frac{\p{\tilde
x^j}}{\p{ x^{i}}}\left(\tilde A_j+\tilde h_{jk}\frac{\p{\tilde
x^k}}{\p{u}}\right)$, hence $\frac{\p{\tilde x^k}}{\p{u}}=0$. We
conclude that $h=\tilde h$. The converse statement is trivial.

\section{Curvature, holonomy and Petrov type}

Let $g$ be an Einstein metric of the form \eqref{Walker} with
$\Lambda\neq 0$, $A=0$ and $H=\Lambda v^2+H_0$, i.e. $H_1=0$.
Consider the vector fields $p=\p_v$, $q=\p_u-\frac{1}{2}H\p_v$ and
the distribution $E$ generated by  $\p_{x}$ and $\p_y$. We will
use the identification $\Lambda^2\Real^{1,3}\simeq\so(1,3)$. In
\cite{GLE} it is shown that the curvature tensor $R$ of the metric
$g$ is given by
$$ R(p,q)=\Lambda p\wedge q ,\quad R(X,Y)=\Lambda X\wedge Y,\quad
R(X,q)=-p\wedge T(X),\quad R(p,X)=0,$$ where  $X,Y$ are sections
of $E$, and $T(X)=-R(X,q)q$ is a symmetric endomorphism of $E$. It
holds $T^j_i=-R^j_{4i4}$. Obviously, the metric $g$ is
indecomposable if and only if $T\neq 0$. In this case the holonomy
algebra at any point $m\in M$ equals to $\Real p_m\wedge
q_m+\so(E_m)+p_m\wedge E_m$ and it coincides with the maximal
subalgebra $\simil(2)\subset \so(1,3)$ preserving the null line
$\Real p_m$. If the $T$ is identically zero, then the manifold is
decomposable and the holonomy algebra coincides with $\Real
p_m\wedge q_m\oplus \so(E_m)$. The holonomy algebras of the
Einstein Walker metrics are found also in each of the papers
\cite{Schell,H-L,GalEinsteinHol}. Remark  that the manifolds under
the consideration never admit (even locally) any null parallel
vector field.

For the Weyl tensor we get $$ W(p,q)=\frac{\Lambda}{3} p\wedge
q,\, W(p,X)=-\frac{2\Lambda}{3} p\wedge X,\,
W(X,Y)=\frac{\Lambda}{3} X\wedge Y,\, W(X,q)=-\frac{2\Lambda}{3}
X\wedge q -p\wedge T(X).$$ In \cite{H-L,HallBook} it is shown that
the Petrov type of the metric $g$ is either II or D (and it may
change
 from  point to point), in particular, the manifolds under the
 consideration are {\it algebraically special} in the sense of the Petrov classification \cite{Petrov}.
Using the Bel criteria and the above expression for $W$, it is
easy to see that $g$ has type II at a point $m\in M$ if and only
if $T_m\neq 0$, and $g$ has type D at a point $m\in M$ if and only
if $T_m=0$. Since $T_m$ is symmetric and trace-free, it is either
zero or it has rank 2. Hence, $T_m=0$ if and only if $\det T_m=0$.

\section{Examples}

{\bf Case $\Lambda<0$.} We consider a partial solution of Equation
\eqref{f,L<0}. This gives an  Einstein metric. Then we find new
coordinates such that $A=0$ and $h$ may depend on~$u$.

\begin{ex} Let $\p_yf=0$, then $f$ satisfies $x^2\p_x^2f-2f=0$. Hence,
$f=\frac{c_1(u)}{x}+x^2c_2(u)$. The function $\frac{c_1(u)}{x}$
defines a Killing form, and we take $f=c(u)x^2$, then $A=2xc(u)\d
y$, and we chose  $H_0=-\Lambda x^4c^2(u)$. To get rid of $A$, we
need to solve the system of equations
$$\frac{d x(u)}{d u}=0,\qquad \frac{d y(u)}{d u}=2\Lambda c( u) x^3( u).$$
Imposing the initial conditions $x(0)=\t x$ and $y(0)=\t y$, we
get that the inverse transformation to the required one has the
form
$$v=\t v,\quad  x=\t x, \quad y=\t y +2\Lambda b(u)\t x^3,\quad u=\t u,$$
where $b(u)$ is the function such that $\frac{ \d b(u)}{\d
u}=c(u)$ and $b(0)=0$.  With respect to the obtained coordinates,
we get $$ g=2\d v\d u+h(u)+\big(\Lambda v^2+3\Lambda
x^4c^2(u)\big)(\d u)^2,$$
$$h(u)=\frac{1}{-\Lambda\cdot x^2}\Big(\big(36\Lambda^2b^2(u)x^4+1\big)(\d x)^2
+12\Lambda b(u)x^2\d x\d y+(\d y)^2\Big).$$ Let $c(u)\equiv 1$,
then $b(u)=u$. It holds $\det T=-9\Lambda^4 x^4\left(x^4+
v^2\right)$. This shows that $\det T_m=0$ $(m=(v,x,y,u))$ if and
only if $v=0$. In particular, the metric is indecomposable.  The
metric $g$ is of Petrov type D on the set $\{(0,x,y,u)\}$ and it
is of type II on its complement.

If $b(u)=u$, then the Lie algebra of Killing vector fields of the
obtained metric is spanned by the vector fields $\p_y$,
$2v\p_v+x\p_x+y\p_y-2u\p_u$, $\p_u-2\Lambda x^3\p_y$.
\end{ex}

\begin{ex} The functions $f=x^2y$  and $H_0=-\Lambda x^4y$ are partial solutions
of \eqref{f,L<0}. Then, $A=-x^2\d x \d u +2xy\d y$. Consider the
transformation with the  inverse one given by
$$v=\t v,\,\, x=\t x(1+3\Lambda \t u\t x^3)^{-\frac{1}{3}},\,\, y=\t y(1+3\Lambda \t u \t x^3)^{\frac{2}{3}},\,\, u=\t u.$$
  With respect to the obtained coordinates, we get
$$g=2\d v\d u+h(u)+
 \Lambda\left( v^2+3x^4y^2+\frac{x^6}{\rho^2}\right)(\d
 u)^2,\qquad \rho=1+3\Lambda ux^3,$$ $$h(u)=\frac{1}{-\Lambda}\left(
\left(36\Lambda^2x^2y^2u^2+\frac{1}{x^2\rho^2}\right)(\d x)^2
+12\Lambda\rho yu\d x\d y+\frac{\rho^2}{x^2}(\d y)^2\right).$$ It
can be checked that $\det T<0$ everywhere, hence the metric $g$ is
indecomposable and it is of Petrov type II everywhere.

The Lie algebra of Killing vector fields of the obtained metric is
spanned by the vector fields  $3v\p_v+x\p_x+y\p_y-3u\p_u$ and
$\Lambda x^4\p_x-2\Lambda x^3y\p_y+\p_u$.
\end{ex}

The above two metrics are not isometric after any change of the
functions $H_0$, see Section \ref{Seccoord}.

{\bf Case $\Lambda>0$.}  Assuming $\p_yf=0$, we get the solution
$$f=c_1(u)\cos x\,  +c_2(u)\left(\ln\left(\tan\frac{x}{2}\right)\cos x\,  +1\right)$$ of Equation
\eqref{fL>0}. Obviously, $yf$ is a solution of Equation
\eqref{fL>0} as well. The function $c_1(u)\cos x$ defines a
Killing form $A$, and we omit it.

\begin{ex}
The function  $f=\ln\left(\tan\frac{x}{2}\right)\cos x\,  +1$ is a
solution of the first equation from \eqref{fL>0}. We get
$A=\left(\cos x\,  -\ln\left(\cot\frac{x}{2}\right)\sin^2
x\,\right)\d y$. Consider the transformation
$$\t v= v,\quad\t x= x,\quad
\t y= y -\Lambda u
\left(\ln\left(\tan\frac{x}{2}\right)-\frac{\cos x}{\sin^2
x}\right),\quad \t u= u.$$
  With respect to the obtained coordinates, we get
$$ g=2\d v\d u+h(u) +\left(\Lambda v^2+\tilde H_0\right)(\d
u)^2,$$ $$h(u)=\left(\frac{1}{\Lambda}+\frac{4\Lambda u^2}{\sin^4
x }\right)(\d x)^2+\frac{4u}{\sin x\, } \d x \d y+\frac{\sin^2
x\,}{\Lambda}(\d y)^2,$$ where $\t H_0$ satisfies $\Delta_h\t
H_0=-\frac{1}{2}h^{ij}\ddot h_{ij}$. An example of such $\t H_0$
is $$\t H_0=-\Lambda\left(\frac{1}{\sin^2
x\,}+\ln^2\left(\cot\frac{x}{2}\right)\right).$$ Coming back to
the initial coordinates, we get
$H_0=\Lambda\cdot\left(\ln\left(\tan\frac{x}{2}\right)\cos x\,
+1\right).$  We have $$\det
T=-\frac{\Lambda^4}{\sin^4x}\left(v^2+\left(\ln\left(\cot\frac{x}{2}\right)\cos
x\,  -1\right)^2\right).$$ Hence $g$ is of Petrov type D on the
set $\left\{(0,x,y,u)|\ln\left(\cot\frac{x}{2}\right)\cos x\,
-1=0\right\}$ and it is of type II on the complement to this set.
The metric is indecomposable.

The Lie algebra of Killing vector fields of the obtained metric is
spanned by the vector fields $\p_y$ and
$\p_u+\Lambda\left(\frac{\cos x\,  }{\sin^2
x\,}-\ln\left(\tan\frac{x}{2}\right)\right)\p_y$.
\end{ex}

\begin{ex} The functions  $f=y\cos x$ and
$H_0=\Lambda\cdot\big(-y^2\cos^2  x\,
+\ln\left(\tan\frac{x}{2}\right)\big)$  are partial solutions of
\eqref{fL>0}. We get $A=-\cot x\d x-y\sin^2 x\,\d y$. Consider the
transformation
$$\t v= v,\quad\t x=\arccos\left(e^{\Lambda u}\cos x\,  \right),\quad\t y= ye^{-\Lambda u},\quad\t u= u.$$
  With respect to the obtained coordinates, we get
$$ g=2\d v\d u+h(u) +\Lambda\left( v^2-y^2 e^{-2\Lambda u}+\frac{1}{2}\ln\frac{1-\rho}{1+\rho}  \right)(\d u)^2,\quad
\rho=e^{-\Lambda u}\cos  x,$$
$$h(u)=\frac{e^{-2\Lambda u}\sin^2 x\,}{\Lambda\cdot\left(1-\rho^2 \right)}(\d x)^2+\frac{e^{2\Lambda u}(1-\rho^2)}{\Lambda}(\d
y)^2.$$
 It holds $\det
T=-\frac{\Lambda^4\cdot \left(4y^2\cos^2 x +\left(\rho
+2v\right)^2\right)}{4\left(1-\rho^2\right)^2}.$ Hence the metric
is indecomposable and $T$ vanishes at a point $(v,x,y,u)$ if and
only if $y\cos x\,  =0$ and $e^{-\Lambda u}\cos x\, +2v=0$.
Consequently, the metric $g$ is of Petrov type D at a point
$(v,x,y,u)$ if and only if $x=\frac{\pi}{2}$ and $v=0$, or $y=0$
and $\cos x\, =-2ve^{\Lambda u}$. It is of Petrov type II at other
points.

The Lie algebra of Killing vector fields of the obtained metric is
spanned by the vector field $-\Lambda \cot x\p_x-\Lambda
y\p_y+\p_u$.
\end{ex}

The metrics of the last two examples are not isometric after any
change of $H_0$.

 {\small {\bf Acknowledgements.}  I am thankful to D.V.~Alekseevsky  for helpful
suggestions and Thomas Leistner for useful discussions. Most of
the calculations are done using Maple~12. I am thankful to Ian
Anderson for tutorials on Maple. The final form of the paper is
due to helpful remarks of Referee.   The work was supported by the
grant CZ.1.07/2.3.00/20.0003. }

%\section*{References}

%\vskip1cm

\end{document}